\documentclass[runningheads]{llncs}
\usepackage[T1]{fontenc}
\usepackage{newtxtext}       %
\usepackage[varvw]{newtxmath}       

\usepackage{xcolor,amsmath,bbm}
\newtheorem{assumption}{Assumption}

\usepackage{algorithm}
\usepackage{algpseudocode}

\def \bE {\mathbb{E}}
\def \bP {\mathbb{P}}
\def \bR {\mathbb{R}}

\def \bd {d}
\def \bx {x}
\def \by {y}
\def \bv {v}
\def \bp {p}

\def \mF {\mathcal{F}}
\def \mT {\mathcal{T}}

\def \bEpiO {\mathbb{E}_{\pi}}
\def \conv {\text{conv}}

\def \proj {\text{proj}}

\def \dist {\text{dist}}
\def \dom {\text{dom}}

\def \ri {\text{ri}}

\def \GoalVio {\textbf{GoalVio}}
\def \Reward {\textbf{Reward}}

\usepackage{graphicx}
\usepackage{booktabs} 

\begin{document}
\title{Online Nonconvex Decision Making with Approximate Oracles}
\titlerunning{Online Nonconvex Decision Making}
%
\author{Rui Chen \and Oktay G\"unl\"uk \and Andrea Lodi \and Guanyi Wang}
%
%
\institute{}
%
\maketitle              
\begin{abstract}
We study online decision making with possibly nonconvex (for example, combinatorial or integer) local feasible sets and convex long-term constraints. At each period, a decision maker selects a local action to maximize cumulative reward while ensuring that the sequence of actions satisfies the global constraints in the long term. It is known that $O(\sqrt{T})$ long-term constraint violation and $O(\sqrt{T})$ reward regret are achievable if given access to a local optimization oracle that exactly solves optimization problems over each local feasible set. Instead of requiring exact optimization over each local feasible set, we assume access to an approximate local optimization oracle that permits optimization errors quantified by two parameters $(\delta_+,\delta_\times)$. We develop an approximate online dual mirror descent algorithm and show that its guarantees depend on a feasibility slack parameter $\underline d$, which measures the compatibility between local feasible sets and the long-term constraint. Under the standard random permutation model, when $\delta_+=O(1)$, $\delta_\times<\underline d$ and $\underline d-\delta_\times=\Omega(1)$, the algorithm achieves $O(\sqrt T)$ long-term constraint violation and reward regret up to $O\bigl(\sqrt T+\delta_\times T\bigr)+\delta_+ T$. We also show a matching limitation: when $\delta_+=0$ and $\delta_\times>\underline d$, no online algorithm relying on such approximate oracle outputs can guarantee sublinear long-term constraint violation in general.

\end{abstract}
\section{Introduction}\label{sec:Intro}
Online platforms increasingly rely on sequential decision-making algorithms to allocate resources, match users, and control exposure in real-time \cite{mehta2010online,devanur2013randomized,li2010contextual,liu2019globally}. Examples include assigning tasks to workers \cite{tong2017flexible}, allocating inventory across regions \cite{asadpour2020online}, and selecting advertisements or assortments for users \cite{devanur2009adwords,gong2022online}. In such settings, platforms make irrevocable decisions over time to optimize immediate objectives, such as revenue or user engagement. Decisions must be made before the future is observed, often with little time to even solve the local problem exactly. At the same time, concerns about the long-term effect of local decisions have become increasingly prominent in algorithmic decision making \cite{calders2009building,dwork2012fairness,patro2022fair}. For example, workload might be allocated unevenly to workers to optimize local efficiency while fairness over time must be taken into consideration \cite{lodi2024fairness,lodi2024framework}.

To address this challenge, we consider the following online decision making model. A decision maker makes a sequence of irrevocable decisions over a finite horizon of \( T \) time steps. At each time step \( t \in \{1,\dots,T\} \), the decision maker observes a local feasible set \( \Omega_t \) and must select a decision
\[
(r_t, x_t, y_t) \in \Omega_t .
\]
Here, \( r_t \in \mathbb{R} \) denotes the immediate reward obtained at time \( t \), \( x_t \) represents auxiliary local decision variables, and \( y_t \in \mathbb{R}^m \) quantifies the \textit{long-term impact} of the decision at the $t$-th time step to the aggregated long-term goals. 

We model our sequential decision-making environment with two types of constraints: We use set \( \Omega_t \) to capture physical or logical constraints that must be satisfied at every time step. We assume the sets \( \Omega_t \) are compact but may be non-convex, for example, they might be mixed integer sets representing certain combinatorial structures, where exact optimization can be time-consuming or even practically infeasible.
In contrast, long-term constraints are not enforced on a per-decision basis but rather on the aggregate {global} outcomes. We assume that \( \Psi \subseteq \mathbb{R}^m \) is a closed convex set representing acceptable long-run outcomes. The decision maker aims to ensure that the aggregated long-term impact satisfies 
\begin{equation}\label{global_constr}
	\frac{1}{T}\sum_{t=1}^T y_t \in \Psi \quad \text{or} \quad \sum_{t=1}^T y_t \in T \Psi ~.    
\end{equation}
Note that we do \textit{not} assume that the long-term goal set $\Psi$ is bounded. This modeling choice generalizes standard resource allocation problems, where \( \Psi \) is typically a polyhedral packing set, by permitting constraints that are nonlinear, nonseparable, and potentially unbounded.

The offline benchmark problem is therefore given by
\begin{equation}
	\begin{aligned}
		z^*=\max_{(r_t,x_t,y_t)} \quad & \sum_{t=1}^T r_t \\
		\text{s.t.} \quad & (r_t,x_t,y_t) \in \Omega_t , \quad t=1,\dots,T, \\
		& \sum_{t=1}^T y_t \in T\Psi .
	\end{aligned}
	\label{offline}
\end{equation}
This formulation captures the trade-off between local rewards and long-term constraints: maximizing $\sum_t r_t$ may require taking actions that drive the average long-term impact $\frac{1}{T}\sum_t y_t$ away from $\Psi$. The decision maker must therefore navigate this trade-off dynamically. 

\section{Problem Setting and Related Work}\label{sec:prob}

\subsection{Performance Metrics}

In the online setting, the decision maker observes \( \Omega_t \) only at time \( t \) and must choose \( (r_t,x_t,y_t) \) without knowledge of future feasible sets. Because future constraints are unknown, it is generally impossible to guarantee that the long-term constraint is satisfied in hindsight in general. Therefore, following the literature \cite{agrawal2014fast,gupta2016experts}, we treat the long-term constraint \eqref{global_constr} as a \textit{soft} constraint and evaluate an online algorithm using two performance metrics.

\paragraph{Reward.}  
The accumulated reward of the online policy is
\[
\Reward := \sum_{t=1}^T r_t .
\]

\paragraph{Goal Violation.}

The deviation from the long-term goals is measured by
\[
\GoalVio := \mathrm{dist}_{\|\cdot\|}\!\left(\sum_{t=1}^T y_t,\, T\Psi \right),
\]
where \( \mathrm{dist}_{\|\cdot\|}(z, A):=\inf_{a\in A}\|z-a\| \) denotes the distance induced by $\|\cdot\|$ to a closed set with $\|\cdot\|$ being a generic norm on the impact space $\bR^m$. 
Let $\|\cdot\|_*$ denote the dual norm of $\|\cdot\|$.

\subsection{Problem Setting}
\paragraph{Motzkin-decomposable goal set $\Psi$.}
To allow more general goals, we assume throughout the paper that $\Psi$ admits a Motzkin decomposition. This is slightly more general than the typical assumption that the goal set $\Psi$ has to be bounded \cite{agrawal2014fast}.
\begin{assumption}
	The goal set $\Psi$ satisfies
	\begin{equation}
		\Psi = Q+C,
	\end{equation}
	where $Q\subseteq \bR^m$ is nonempty, compact, and convex, and $C\subseteq \bR^m$ is a closed convex cone.
\end{assumption}
This structure allows unbounded long-term goal sets while keeping certain operations (in particular, linear optimization on certain directions) over $\Psi$ well-defined. Let $C^\circ$ denote the polar cone of the recession cone $C$ of $\Psi$, i.e.,
$$C^\circ:=\{p\in\bR^m:\langle p,c\rangle\leq 0~\forall c\in C\}.$$

\paragraph{Regularity assumptions.} To derive quantitative guarantees on long-term goal violation and reward, we further impose the following regularity assumptions throughout the paper. In particular, we consider a feasibility slack parameter $\underline{d}$, which measures how much one can perturb a convex combination of local decisions so that it remains acceptable with respect to the long-term constraint. This condition does not require any individual local action to be strictly feasible, but it is stronger than feasibility of the aggregate offline problem. The margin provides room for the algorithm to compensate for errors in local optimization.
\begin{assumption}\label{assmpt:regularity}
	There exist constants $$d_y=O(1),~d_r=O(1),~\text{and}~\underline{d}=\Omega(1)$$ such that the following hold for all time steps $t=1,\ldots,T$:
	\begin{enumerate}
		\item \textit{Uniform boundedness for local decisions.} For every decision $(r,x,y)\in\Omega_t$,
		\begin{align*}
			\max_{v\in Q}\|y-v\|\leq d_y,~r\in[0,d_r].
		\end{align*}
		\item \textit{Strong feasibility in the norm $\|\cdot\|$.} There exists $(\tilde{r}_t,\tilde{\bx}_t,\tilde{\by}_t)\in \conv(\Omega_t)$ such that $$\{\tilde{\by}_t\}+\underline{d}B_{\|\cdot\|}\subseteq \Psi,$$ where $B_{\|\cdot\|}$ denotes the unit ball in the norm $\|\cdot\|$. 
	\end{enumerate}
\end{assumption}

\paragraph{Approximate local optimization oracles.} A central difficulty in our setting is that the local feasible set faced in each period may be non-convex and hard to optimize over. In such cases when the exact oracle (for example, via a mixed-integer solver) may be computationally expensive, and practitioners may have to settle for an approximate optimization algorithm or a suboptimal solution with a certified but nonzero optimality gap. To model such a compromise, we assume the following approximate local optimization oracles.
\begin{assumption}\label{assmpt:approx}
	There exist constants $\delta_+\geq 0$ and $\delta_\times\geq 0$ such that, at  each time step $t$, we have one access to an approximate oracle $O_t(\cdot;\Omega_t)$ that inputs any vector $p_t\in C^\circ$, and outputs an approximate solution $(r_t,x_t,y_t)\in\Omega_t$ that is $(\delta_+,\delta_\times)$-approximate with respect to the $p_t$-priced subproblem. Specifically, denoting
    \begin{equation}\label{define_et}
        e_t:=\max_{(r,x,y)\in\Omega_t} \{r-\langle p_t, y\rangle\}-[r_t-\langle p_t, y_t\rangle],
    \end{equation}
    we assume
    \begin{equation}\label{equation:local_approx}
		e_t\leq \delta_++\delta_{\times}\|p_t\|_*.
	\end{equation}
\end{assumption}
\begin{remark}
    Assumption \ref{assmpt:approx} admits oracles with bounded additive optimization errors (quantified by $\delta_+$). Moreover, assume $0\in Q$, then one can easily show that Assumption \ref{assmpt:approx} admits oracles with multiplicative optimization errors under the regularity assumptions imposed by Assumption \ref{assmpt:regularity}. For example, assuming a nonnegative(/nonpositive) optimal objective, outputs of oracles with a $1-(/+)\epsilon$ approximation ratio to the $p_t$-priced subproblem\begin{align}\label{pt-pricedsub}
        \max_{(r,x,y)\in\Omega_t} \{r-\langle p_t, y\rangle\}
    \end{align}
    are $(\delta_+,\delta_\times)$-approximate with respect to \eqref{pt-pricedsub} with $\delta_+=\epsilon d_r$ and $\delta_\times=\epsilon d_y$.
\end{remark}

\paragraph{Reference function.} To develop guarantees in terms of a generic norm $\|\cdot\|$, we need to consider a proper reference function as follows.
\begin{assumption}\label{assmpt:reference_function}
	Let $\Phi:\bR^m\to\mathbb R$ be differentiable and $1$-strongly convex and $L$-smooth for some $L=O(1)$ with respect to $\|\cdot\|_*$, i.e., for all $p,q\in \bR^m$,
	$$
	\Phi(p)+\langle\nabla\Phi(p),q-p\rangle
	+
	\frac{L}{2}\|q-p\|_*^2
	\geq
	\Phi(q)
	\ge
	\Phi(p)+\langle\nabla\Phi(p),q-p\rangle
	+
	\frac12\|q-p\|_*^2.
	$$
\end{assumption}
For all $p,q\in \bR^m$, define the Bregman divergence $D_\Phi(q,p)$ generated by $\Phi$:
$$
D_\Phi(q,p)
:=
\Phi(q)-\Phi(p)-\langle\nabla\Phi(p),q-p\rangle.
$$
The following result follows directly from Assumption \ref{assmpt:reference_function}.
\begin{proposition}\label{prop:D_norm}
	For all $p,q\in \bR^m$, $D_\Phi(q,p)=\Theta(\|q-p\|_*^2)$, and specifically,
	$$\frac{L}{2}\|q-p\|_*^2\geq D_{\Phi}(q,p)\geq\frac{1}{2}\|q-p\|_*^2.$$
\end{proposition}

\subsection{Related Work}
Early studies of online decision making problems include online bipartite matching \cite{karp1990optimal}, routing \cite{awerbuch1993throughput}, secretary problems \cite{dynkin1963optimum,kleinberg2005multiple}, advertising \cite{mehta2007adwords} and knapsack \cite{babaioff2007knapsack} problems. Much of the early literature focuses on worst-case competitive analysis. More recently, attention has shifted toward less pessimistic stochastic input models \cite{gupta_singla_2021}. For example, in the context of online bipartite matching \cite{goel2008online}, it is shown that the greedy algorithm achieves a strictly better competitive ratio under the random permutation model than is possible in the worst-case setting.

One of the most well-studied online decision making problems is (multiple-choice packing) \textit{Online Linear Programming} (OLP), also known as online resource allocation. It is a special case of \eqref{offline} with a polyhedral goal set $\Psi=\{\by:\by\leq\bd\}$ where $\bd>\mathbf{0}$ represents a vector of resource budgets. At each time step $t$, the local feasible set $\Omega_t$ takes the form
\begin{equation}
	\Omega_t=\Big\{(r,\bx,\by):r={\alpha}_t^\top\bx,
	\bx\in \Delta,~\by={A}_t\bx\Big\},\label{eq:OLP}
\end{equation}
where ${A}_t$ and $\alpha_t$ are nonnegative, and $\Delta$ is the standard simplex $\{\bx\geq {\bf 0}:\sum_j x_j\leq 1\}$ representing a multiple-choice setting with the option of a void decision $\bx=\mathbf{0}$. 
A stream of work on OLP develops reward guarantees under stochastic input models while maintaining exact feasibility with respect to resource budgets \cite{feldman2010online,molinaro2014geometry,kesselheim2014primal,gupta2016experts,li2022simple}. These results exploit the structure of packing constraints, whereas our setting permits more general convex constraints on aggregate outcomes. In a related line of work \cite{balseiro2023best}, a mirror-descent-based primal-dual algorithm is proposed for a class of online resource allocation problems with long-term packing constraints, allowing nonconvex local decisions.

Online linear programming is substantially generalized in the seminal work on \textit{online stochastic convex programming} \cite{agrawal2014fast}, which allows general bounded local feasible sets \(\Omega_t\) and a bounded convex goal set \(\Psi\). 
Under this framework, a primal-dual online algorithm achieves sublinear reward regret and long-term constraint violation in expectation under stochastic input models. This framework broadens the range of long-term goals that can be incorporated into online decision making, but the feasibility safeguards available for resource budgets do not directly extend to general convex goals.

Our work bridges these lines of research by studying approximate local optimization under more general convex long-term goals, including a class of unbounded goal sets. We establish reward and constraint-violation guarantees and characterize when local approximation errors remain compatible with long-term feasibility. In particular, our results show that the reliability of approximate local decisions depends not only on oracle accuracy but also on the feasibility slack available in the problem.

\section{A Primal-Dual Framework}\label{sec:alg}
In this section, we develop a primal-dual framework for the online decision-making problem introduced in Section 2. The key idea is to guide local decisions via dual multipliers associated with the long-term constraints, whose values are learned online. This perspective allows us to coordinate non-convex, irrevocable actions over time while controlling long-term constraint violations.


\subsection{Convexification and a Dual Interpretation}\label{subsec:convexification}
The offline problem in \eqref{offline} is generally non-convex due to the local feasible sets \( \Omega_t \). To facilitate our algorithm development and analysis, we consider a partial convexification of \eqref{offline} obtained by replacing each \( \Omega_t \) with its convex hull \( \operatorname{conv}(\Omega_t) \), namely,
\begin{equation}
	\begin{aligned}
		z_R := \max_{(r_t,x_t,y_t)} \quad 
		& \sum_{t=1}^T r_t \\
		\text{s.t.} \quad 
		& (r_t,x_t,y_t) \in \operatorname{conv}(\Omega_t), \quad t=1,\dots,T, \\
		& \sum_{t=1}^T y_t \in T\Psi .
	\end{aligned}
	\label{offline_convex}
\end{equation}
Problem \eqref{offline_convex} is a convex relaxation that serves as a natural benchmark for online algorithms. By construction, \(z_R \ge z^*\), where \(z^*\) denotes the optimal value of \eqref{offline}; moreover, \(z_R=z^*\) when all \(\Omega_t\) are convex. Importantly, this relaxation preserves the structure of the long-term constraints. We assume throughout that problem \eqref{offline_convex} is feasible. Finally, while non-convexity of \((\Omega_t)_{t=1}^T\) can create a duality gap between \(z_R\) and \(z^*\), the relative gap (normalized by $1/T$) typically shrinks as \(T\) grows under mild conditions \cite{aubin1976estimates,dey2026asymptotically}, as a consequence of the Shapley-Folkman lemma \cite{starr1969quasi}.

For each time step \( t \), define the concave function
\[
f_t(y) := \max \{ r : (r,x,y) \in \operatorname{conv}(\Omega_t) \},
\]
which captures the maximum achievable reward at time \( t \) given a contribution \( y \) to the long-term constraints. Using this notation, problem \eqref{offline_convex} can be written compactly as
\begin{equation}\label{offline_convex_primal}
	z_R=\max_{y_1,\dots,y_T}~
	\sum_{t=1}^T f_t(y_t) -\mathbb{I}_{T\Psi}\left(\sum_{t=1}^T y_t\right),
\end{equation}
where $\mathbb{I}_{{T\Psi}}$ is the indicator function of $T\Psi$ defined as 
\begin{displaymath}
	\mathbb{I}_{{T\Psi}}(\by)=\begin{cases} 0&\text{if}~\by\in {T\Psi},\\+\infty&\text{otherwise.}\end{cases}
\end{displaymath}

The formulation in \eqref{offline_convex_primal} highlights a natural tradeoff between reward and constraint violation. Introducing a dual vector \( p \in \mathbb{R}^m \) associated with the aggregated fairness impact $\sum_{t=1}^Ty_t$ and considering the Fenchel dual \cite{rockafellar1970convex} of \eqref{offline_convex_primal}, we obtain a dual problem of the form
\begin{equation}\label{offline_convex_dual}
	z_D=\min_{p} \;
	T h_\Psi(p) - \sum_{t=1}^T f_t^*(p),
\end{equation}
where \[ h_\Psi(p) = \sup_{v \in \Psi} \langle p, v\rangle \] is the support function of \( \Psi \), and \( f_t^* \) denotes the conjugate of \( f_t \), i.e.,
$$f_t^*(p)=\min_{(r,x,y) \in \conv(\Omega_t)} \langle p, y\rangle -r=\min_{(r,x,y) \in \Omega_t} \langle p, y\rangle-r.$$ Importantly, evaluating \( f_t^*(p) \) requires exactly optimizing a linear function $\langle p, y\rangle-r$ over $\Omega_t$. In our algorithmic development, we are going to replace this by an approximate solution obtained from an oracle call to the local optimization oracle. 


\subsection{An Approximate Dual Mirror Descent Algorithm}\label{subsec:algorithm}
We now present an approximate online dual mirror descent algorithm for problem \eqref{offline} introduced in Section \ref{sec:prob}. The algorithm is motivated by the dual interpretation in Section \ref{subsec:convexification}. From a high-level perspective, the algorithm maintains a vector of dual multipliers \( p_t \in \mathbb{R}^m \) to penalize long-term constraint violations. The local decisions are generated from a call to the approximate local optimization oracle on a $p_t$-perturbed local problem. The dual multipliers are updated with an approximate mirror descent step.

At each time step \( t \), after observing the local feasible set \( \Omega_t \), given the current dual multipliers $p_t$, the algorithm selects a decision by calling the oracle $O_t(p_t;\Omega_t)$ that approximately solves the local optimization problem
\begin{equation}\label{step:local}
	\max_{(r,x,y)\in\Omega_t}
	\bigl\{ r - \langle p_t, y \rangle \bigr\}.
\end{equation}
We denote an optimal solution of \eqref{step:local} by $(\bar r_t,\bar x_t,\bar y_t)$, and the approximate solution obtained from the oracle output by $(r_t, x_t, y_t)$.

The dual multipliers are then updated. Specifically, for $p_t\in C^{\circ}$, let
\begin{equation}\label{step:extreme}
	v_t \in \arg\max_{v \in Q} \langle p_t, v\rangle\subseteq \arg\max_{v \in \Psi} \langle p_t, v\rangle
\end{equation}
be a maximizer of the long-term goal set in the direction \( p_t \). The dual update is given by a mirror descent-like step\begin{equation}\label{step:mirror_descent}
	p_{t+1}\in\arg\min_{p\in C^\circ}\{\eta_t \langle p,v_t-y_t\rangle+D_\Phi(p,p_t) \},
\end{equation}
where \( \eta_t > 0 \) is a stepsize. Intuitively, if the realized impact \( y_t \) exceeds the fairness \( v_t \) along some dimension, the corresponding dual price increases, discouraging similar decisions in the future; under-utilized dimensions receive lower prices, incentivizing compensating actions.
Algorithm 1 summarizes the procedure. {Note that Algorithm \ref{alg:LagDual} does not need the length of the horizon $T$ or the oracle-error parameters as inputs.}

\begin{algorithm}[tb!]
	\caption{An Approximate Online Dual Mirror Descent Algorithm for Problem \eqref{offline}}
	\label{alg:LagDual}
	\begin{algorithmic}
		\State {Initialize} $\bp_1={0}$\;
		\For{$t=1,\ldots,T$}
		\State Select the local decision via the approximate optimization oracle:
		$$(r_t,x_t,y_t)\leftarrow O_t(p_t;\Omega_t)$$
		\State Set $v_t\in\arg\max_{v\in Q}  \langle p_t, v\rangle$ 
		\State Update the dual multipliers via approximate mirror descent:
		$$p_{t+1}\in\arg\min_{p\in C^\circ}\{\eta_t \langle p,v_t-y_t\rangle+D_\Phi(p,p_t) \}$$
		\EndFor
	\end{algorithmic}
\end{algorithm}

Indeed, Algorithm \ref{alg:LagDual} can be interpreted as an instance of approximate online mirror descent applied to the Fenchel dual problem \eqref{offline_convex_dual}, i.e.,
$$\min_p \sum_{t=1}^T \ell_t(p),$$ in which $\ell_t(\cdot) := h_{\Psi}(\cdot) - f_t^*(\cdot)$ is a convex function.
Note that a subgradient of \(\ell_t(\cdot)\) at \(p_t\) is given by $v_t-\bar{y}_t$. As the exact optimal $\bar y_t$ is not available to the algorithm, the dual update takes the form
\[
\arg\min_{p\in C^\circ}\{\eta_t \langle p,\hat{\nabla}\ell_t(p_t)\rangle+D_\Phi(p,p_t) \},
\]
where \(\hat{\nabla} \ell_t(p_t):=v_t-y_t\) is an approximation of the subgradient of the convex function  
\(\ell_t(\cdot)\) evaluated at \(p_t\).
Thus, the algorithm performs approximate online mirror descent on the dual objective \eqref{offline_convex_dual}. Since the dual update is not an exact online mirror descent step \cite{fang2022online}, standard results for online mirror descent do not apply here. Moreover, even if we apply an exact online mirror descent step, the dual guarantees (e.g., regret bounds in terms of the dual objective \eqref{offline_convex_dual}) alone do not directly translate into primal performance guarantees, measured by both the realized primal reward and the constraint violations, with respect to the primal decisions \((r_t,x_t,y_t)_{t=1}^T\). Establishing such primal performance guarantees is therefore the focus of Section \ref{sec:performance-guarantees}.

\section{Primal Performance Analysis} \label{sec:performance-guarantees}
\subsection{Guarantees on Long-Term Constraint Violation}\label{subsec:fairness}

We first show that a strong duality result between the convexified offline problem \eqref{offline_convex_primal} and its dual \eqref{offline_convex_dual} holds.
\begin{lemma}\label{lem:Lag}
	There exist primal and dual solutions $$({r}^*_t,{x}^*_t,{y}^*_t)_{t=1}^{T}\in\prod_{t=1}^{T}\conv(\Omega_t),~p^*\in\bR^m$$ such that $({y}_t^*)_{t=1}^{T}$ optimizes \eqref{offline_convex_primal} and $p^*$ optimizes \eqref{offline_convex_dual} with $z_R=z_D$.
\end{lemma}

Next, we show that Algorithm 1 yields bounded dual multipliers $(p_t)_{t=1}^T$ under the assumptions that $\delta_+=O(1)$ and $\delta_\times<\underline{d}-\Omega(1)$.
\begin{lemma}\label{lem:boundedNorm}
	Suppose $\delta_+=O(1)$, $\delta_\times<\underline{d}$ and $\underline{d}-\delta_\times=\Omega(1)$. If $0<\eta_t\leq\bar{\eta}$ for some $\bar\eta=O(1)$ for all $t$, then Algorithm \ref{alg:LagDual} produces $(p_t)_{t=1}^{T+1}$ with
	$$\max_{t}\|p_t\|_*=O(1).$$
\end{lemma}

We next establish theoretical guarantees on long-term constraint violation achieved by Algorithm \ref{alg:LagDual}. Importantly, these guarantees hold without any assumptions on the input sequence and therefore apply to arbitrary (possibly adversarial) arrival orders of the local feasible sets \((\Omega_t)_{t=1}^T\). Recall that long-term constraint violation is evaluated as a soft long-term constraint on the aggregate outcome, measured by
\[
\GoalVio := \mathrm{dist}_{\|\cdot\|}\!\left(\sum_{t=1}^T y_t,\, T\Psi \right).
\]

\begin{theorem}\label{thm:boundedVio}
    Suppose $\delta_+=O(1)$, $\delta_\times<\underline{d}$ and $\underline{d}-\delta_\times=\Omega(1)$. If $\eta_t=\bar{\eta}/\sqrt{t}$ for some fixed $0<\bar\eta=O(1)$ for all $t$, then Algorithm \ref{alg:LagDual} achieves $$\GoalVio\leq O(\sqrt{T}).$$
\end{theorem}

With $\rho:=\delta_\times/\underline{d}$, Theorem \ref{thm:boundedVio} shows that whenever $\rho<1$, Algorithm \ref{alg:LagDual} may achieve a sublinear long-term constraint violation with proper stepsizes. We next show a one-dimensional example that, if $\rho >1$ and decisions are chosen by $O(p_t;\Omega_t)$ for some sequence $(p_t)_{t=1}^{T}$, then a linear long-term constraint violation cannot be avoided (in the worst case), no matter how $(p_t)_{t=1}^{T}$ are chosen.
\begin{example} \label{exmp:goalvio}
Let $\Psi=(-\infty,1]=Q+C$ with $Q=\{1\}$ and $C=\bR_-$, and $\|\cdot\|=|\cdot|$. Let $\rho>1$ be a fixed constant and $\Omega_t\equiv\bar\Omega$ that consists of two zero-reward actions
$a_0=(r,x,y)=(0,0,0)$ and $a_1=(r,x,y)=(0,0,\rho)$.
One can verify that $\Psi$ and $(\Omega_t)_{t=1}^T$ satisfy Assumption \ref{assmpt:regularity} with $d_y=\max\{1,\rho-1\}$, $d_r=0$, and $\underline{d}=1$. Consider the deterministic oracle that always outputs $O(p;\bar\Omega)=a_1$ for all $p$. Then one can verify that the oracle satisfies Assumption \ref{assmpt:approx} with $\delta_+=0$ and $\delta_\times=\rho$ as for all $p$
$$0-p\cdot\rho\geq\max\{0,0-p\cdot\rho\}-\rho|p|.$$
Therefore, any algorithm relying on such a deterministic oracle to produce the implemented actions would achieve a deterministic linear long-term constraint violation
$\GoalVio=(\rho-1)T=\Omega(T)$, independent of how $(p_t)_{t=1}^T$ are chosen by the algorithm. \\

\noindent
Together with Theorem \ref{thm:boundedVio}, this example establishes a worst-case phase transition in \GoalVio: an $O(\sqrt{T})$ guarantee holds when $\underline d-\delta_\times=\Omega(1)$, whereas $\Theta(T)$ violation can be unavoidable when $\delta_\times>\underline d$.
\end{example}

\subsection{Reward Performance under Uniform Random Permutation}\label{subsec:random_perm}
We now analyze the reward performance of Algorithm \ref{alg:LagDual} under stochastic input models. While analysis in Section \ref{subsec:fairness} guarantees sublinear violation of long-term constraints for any arrival order under proper assumptions on the quality of the local optimization oracles, obtaining meaningful guarantees on reward requires additional assumptions on how the local feasible sets arrive over time. Following standard practice in the online resource allocation literature, we adopt the random permutation model.

Specifically, we assume that there exist \(T\) deterministic feasible sets \(Z_1,\ldots,Z_T\), which may be chosen adversarially and are not known to the algorithm in advance. The decision maker observes these sets in a permuted order: at time step \(t\), the local feasible set is given by
\[
\Omega_t = Z_{\pi(t)},
\]
where \(\pi\) is an unknown permutation of \(\{1,\ldots,T\}\). Owing to the symmetry of the offline problem \eqref{offline}, the optimal offline reward is invariant to the choice of \(\pi\).

The \textit{uniform random permutation model}, also known as the random-order model, assumes that \(\pi\) is drawn uniformly at random from the set of all \(T!\) permutations. This model has been widely studied in online linear programming and resource allocation, and is strictly more general than the IID model, in which each \(\Omega_t\) is sampled independently from a fixed (and possibly unknown) distribution \cite{molinaro2014geometry}.

Throughout this section, we assume Assumption \ref{assmpt:approx} holds for any arrival order. We let \(\mathbb{P}_{\pi}\) and \(\mathbb{E}_{\pi}\) denote the probability measure and the expectation with respect to the random permutation \(\pi\). As a performance benchmark, we consider the optimal objective value \(z_R\) of the partially convexified offline problem \eqref{offline_convex}. By definition, \(z_R \ge z^*\), where \(z^*\) denotes the optimal value of the original offline problem \eqref{offline}. Consequently, any lower bound on the achieved \(\Reward\) relative to \(z_R\) immediately implies the same bound relative to \(z^*\).

We now establish an expected reward guarantee for Algorithm \ref{alg:LagDual} under the uniform random permutation model. 

\begin{theorem}\label{thm:main}
Suppose $\delta_+=O(1)$, $\delta_\times<\underline{d}$ and $\underline{d}-\delta_\times=\Omega(1)$. Let $\eta_t=\bar{\eta}/\sqrt{t}$ for some fixed $0<\bar\eta=O(1)$ for all $t$.
Under the uniform random permutation model, Algorithm \ref{alg:LagDual} achieves
\[
\bEpiO[\Reward]
\;\ge\;
z_R -\delta_+ T-O(\delta_\times T)- O_{\|\cdot\|,m}\!\left(\sqrt{T}\right),
\]
where $O_{\|\cdot\|,m}$ omits some constants that depend on the norm $\|\cdot\|$ and the dimension $m$ of the goal set $\Psi$.
\end{theorem}

\section{Numerical Experiments}\label{sec:numerical}

This section evaluates the numerical performance of the proposed online decision-making algorithm (Algorithm~\ref{alg:LagDual}). We test on an online version of the assortment optimization problem with visibility constraints \cite{barre2025assortment}. Specifically, we examine the following empirical questions: 

\noindent
(1) To what extent does the proposed algorithm outperform or underperform  some baselines such as a greedy algorithm or the offline optimal solution in \Reward ? 

\noindent
(2) How does local oracle accuracy affect long-term constraint violation \GoalVio? How do the observed trends relate to the worst-case results in Section \ref{sec:performance-guarantees}?

The main text presents results on randomly generated instances. Additional results on a deliberately constructed instance, designed to illustrate the worst-case behavior of constraint violation, are deferred to Appendix~\ref{apd:artificial-instance}. 


\noindent 
\textbf{Problem formulation and instance generation.} The objective of assortment optimization is to determine the optimal assortment of products to offer to each customer in order to maximize the total expected revenue. We consider $m$ products in total and $T$ customers arriving online. The corresponding offline problem is NP-hard and takes the form\begin{equation}\label{eq:num-problem}
 \max_{S_1,\ldots,S_T\in\mathcal S}\ \sum_{t=1}^T r_t(S_t)
 \quad\text{s.t.}\quad
 \sum_{t=1}^T \mathbf{1}(i\in S_t)\ge Tf_i,\quad i=1,\ldots,m,
\end{equation}
where
\[
\mathcal S=\{S\subseteq\{1,\ldots,m\}:|S|\le s\}
\]
Here $S_t$ denotes the selected assortment displayed to customer $t$, $s$ is the maximum assortment size, and $f_i$ is the minimum required display frequency of product $i$.

We assume that the local reward function $r_t(\cdot)$ of customer $t$ is determined by a Multinomial Logit model \cite{luce2005individual} whose preference weights are randomly generated. In the notation of our framework, 
\[
\Omega_t=\{(r_t(S),\mathbf1_S,\mathbf1_S):S\in\mathcal{S}\},\quad \Psi=\{f\}+\mathbb R_+^m,
\]
where $f=(f_1,\ldots,f_m)$ and $\mathbf1_S$ denotes the incidence vector of $S$. We use the $\ell_\infty$-norm $\|\cdot\|=\|\cdot\|_\infty$ to measure long-term constraint violation, so that
\[
\GoalVio=\max_{i=1,\ldots,m}\left[Tf_i-\sum_{t=1}^T\mathbf{1}(i\in S_t)\right]_+.
\]
The feasibility slack is then
\[\underline{d}=\min\left\{1-\max_i f_i,\frac{s-\sum_{i=1}^m f_i}{m}\right\}.\]
We generate five base instances with $m=8$, $s=4$ and $T=40000$, and truncate the customers at horizons $T\in\{500,1000,2000,4000,8000,16000,24000,40000\}$ to produce instances with different scales. Further details of instance generation are presented in Appendix \ref{apd:details-random-instances}. We choose relatively small $m$ and $s$ so that the offline problem is solvable. These experiments are intended to illustrate the algorithm’s behavior, rather than to establish computational scalability in the number of products.

\noindent\textbf{Algorithms and implementation.} We compare Algorithm \ref{alg:LagDual} with two benchmarks:\begin{itemize}
    \item A \textit{greedy} algorithm that, at each time step, maximizes the current revenue among assortments permitting completion of the remaining visibility quotas.
    \item The \textit{offline} benchmark which uses the entire customer sequence to compute the hindsight optimal solution.
\end{itemize}
For Algorithm \ref{alg:LagDual}, we use $\Phi(p)=\frac m2\|p\|_2^2$ and $\eta_t=\frac{0.05m}{\sqrt t}$. The local problem \eqref{step:local} reduces to a single-customer assortment optimization problem under a mixture of independent demand and multinomial logit models with a cardinality constraint, which can be approximately solved to any accuracy using a fully polynomial-time approximation scheme (FPTAS) developed in \cite{cao2023revenue}. The oracle settings are chosen to satisfy Assumption \ref{assmpt:approx} with $\delta_+=0.05$ and different values of $\delta_\times$. We summarize oracle accuracy relative to feasibility slack by $\alpha:=\delta_{\times}/\underline{d}$.\\

\noindent
\textbf{Numerical results and discussions.} Figure 1 reports reward regret relative to the offline optimum,
\[
\operatorname{Regret}_R(T)
=
z^*(T)-\mathbf{Reward}(T),
\]
for Algorithm \ref{alg:LagDual} (“Alg. 1”) and the greedy baseline (“Greedy”) across horizons and oracle settings. Results are averaged over the five base instances.

Across the tested settings, Algorithm \ref{alg:LagDual} achieves substantially lower mean reward regret than the greedy baseline. The right-hand panel also indicates that less accurate local optimization, corresponding to larger \(\alpha\), generally increases reward regret. This trend is qualitatively consistent with Theorem 2.

Meanwhile, we observe that the mean constraint violation of Algorithm 1 remains small (ranging from approximately 7 to 14), with no clear monotonic dependence on either \(T\) or \(\alpha\). Thus, on this instance family, the algorithm achieves small long-term constraint violation even under relatively loose oracle settings.
Indeed, the worst-case failure in Example \ref{exmp:goalvio} need not occur for the FPTAS on randomly generated instances. Appendix C.2 therefore examines a deliberately constructed instance to illustrate the dependence of constraint violation on the horizon and the approximation-to-slack ratio.



\begin{figure}[tbp]
	\centering
	\includegraphics[width=\textwidth]{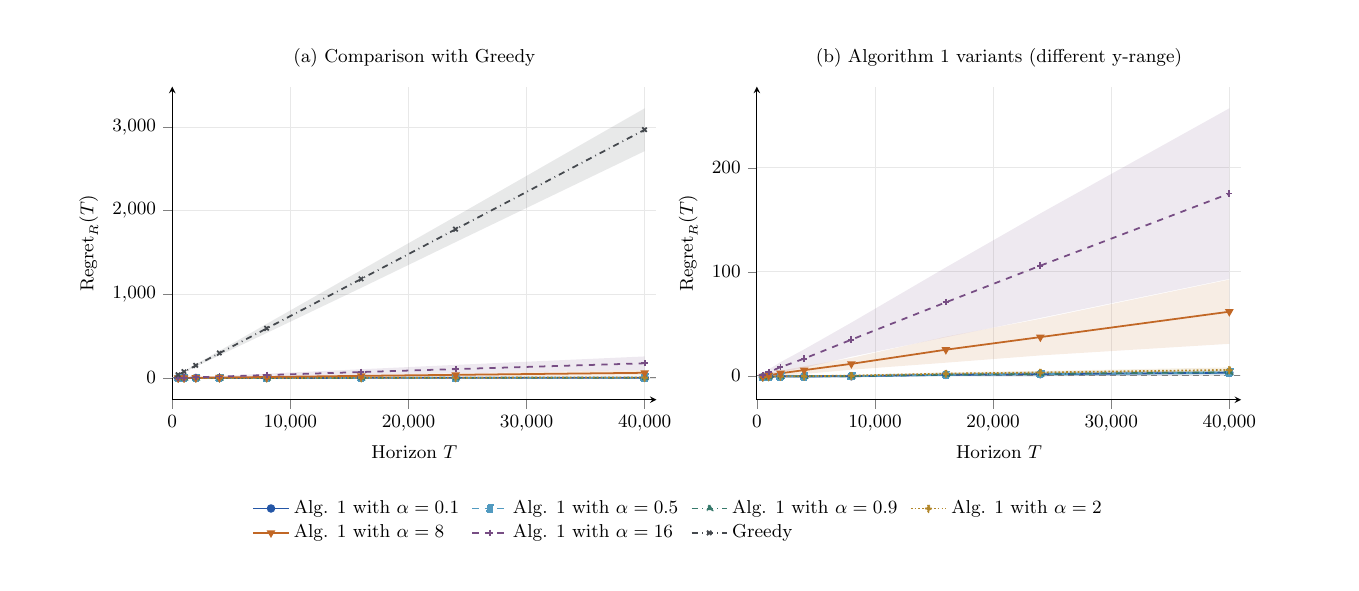}

    \vspace{- 4 ex} 
    \caption{\small Regret versus horizon $T$ over randomly generated instances. Means over five independent seeds, shading denotes pointwise 95\% Student-$t$ intervals. $\mathrm{Regret}_R(T)=z^*(T)-\mathrm{Reward}(T)$ uses a shared offline benchmark. }
	\label{fig:regret-vs-T-random}
\end{figure}

\subsubsection{\ackname} 
Large language models were used to refine the writing and assist in extending several results from the Euclidean norm to general norms. The authors independently verified all mathematical arguments and take full responsibility for the paper’s content.

%
%
%
\bibliographystyle{splncs04}
\bibliography{ref}
%




\newpage
\appendix
\section{Some Supporting Inequalities}
We derive some inequalities that are going to be helpful for our analysis.
\begin{lemma}\label{lem:ineq}
	Let \((\bp_t,\bv_t,\by_t)_{t=1}^T\) be generated by Algorithm~\ref{alg:LagDual}.
	Then, for \(t = 1,\ldots,T\), the following inequalities hold:
	\begin{enumerate}
		\item $\|\bp\|_*\,\underline{d}
		\;\le\;
		h_{\Psi}(\bp) - \langle \bp, \tilde{\by}_t\rangle
		$ for \(\bp \in C^\circ\);
		\item $\langle \bp_t - \bp,\bv_t - \by_t\rangle
		\;\le\;
		\frac{D_\Phi (p,p_t)-D_\Phi(p,p_{t+1})}{\eta_t}
		\;+\;
		\frac{\eta_t}{2}\|\bv_t - \by_t\|^2
		$
		for $p\in C^\circ$;
		\item $\langle p_t,v_t-y_t \rangle\geq (\underline d-\delta_\times)\|p_t\|_*-(d_r+\delta_+)$;
		\item $\|p_{t+1}-p_t\|_*\leq\eta_td_y$.
	\end{enumerate}
\end{lemma}
\begin{proof}
	By the strong feasibility assumption in Assumption \ref{assmpt:regularity},
	$$h_\Psi(p)
	=
	\sup_{v\in\Psi}~\langle p,v\rangle
	\geq
	\sup_{v\in\{\tilde{y}_t\}+\underline{d}B_{\|\cdot\|}}\langle p,v\rangle
	=
	\sup_{\|z\|\le 1}
	\langle p,\tilde y_t+\underline{d}z\rangle
	=\langle \bp, \tilde{\by}_t\rangle+\|\bp\|_*\,\underline{d},$$
	which proves the first inequality.
	
	The first-order optimality condition of \eqref{step:mirror_descent} implies the existence of $\xi_t\in N_{C^\circ}(p_{t+1})$ such that
	\begin{equation}\label{eq:optimality}
		\eta_t(v_t-y_t)+\nabla\Phi(p_{t+1})-\nabla\Phi(p_{t})+\xi_t=0.
	\end{equation}
	By the three-point identity of Bregman divergences,\begin{equation}
		D_{\Phi}(p,p_t)-D_{\Phi}(p,p_{t+1})=D_{\Phi}(p_{t+1},p_{t})-\langle p-p_{t+1},\nabla\Phi(p_t)-\nabla\Phi(p_{t+1}) \rangle.
	\end{equation}
	Note that, by definition of the normal cone,\begin{equation}\label{ineq:normal}
		\langle p-p_{t+1},\xi_t\rangle\leq 0.
	\end{equation}
	Combining \eqref{eq:optimality}-\eqref{ineq:normal},
	\begin{multline*}
		D_{\Phi}(p,p_t)-D_{\Phi}(p,p_{t+1})\geq D_{\Phi}(p_{t+1},p_{t})-\eta_t\langle p-p_{t+1}, v_t-y_t\rangle\\
		\geq \frac{1}{2}\|p_{t+1}-p_{t}\|_*^2-\eta_t\langle p-p_{t+1}, v_t-y_t\rangle\\
		=\eta_t\langle p_t-p, v_t-y_t\rangle+\frac{1}{2}\|p_{t+1}-p_{t}\|_*^2-\eta_t\langle p_t-p_{t+1},v_t-y_t\rangle\\
		\geq \eta_t\langle p_t-p, v_t-y_t\rangle-\frac{\eta_t^2}{2}\|v_t-y_t\|^2,
	\end{multline*}
	where the last inequality is due to the Cauchy-Schwarz inequality. This proves the second inequality.
	
	By the uniform boundedness assumption in Assumption \ref{assmpt:regularity} and Assumption \ref{assmpt:approx},
	\begin{equation}\label{ineq:apply_approx}
		\tilde{r}_t+d_r-\langle p_t, y_t\rangle\geq r_t-\langle p_t, y_t\rangle\geq \tilde r_t-\langle p_t, \tilde{y}_t\rangle-\delta_+-\delta_{\times}\|p_t\|_*.
	\end{equation}
	Note that $h_{\Psi}(p_t)=\langle p_t,v_t\rangle$. Applying the first inequality in this lemma with $p=p_t$ to inequality \eqref{ineq:apply_approx} yields the third inequality.
	
	Optimality condition \eqref{eq:optimality} and inequality \eqref{ineq:normal} with $p=p_t$ imply
	\begin{equation}\label{ineq:grad_p}
		\langle\nabla\Phi(p_{t+1})-\nabla\Phi(p_{t}),p_{t+1}-p_t \rangle\leq-\eta_t\langle v_t-y_t,p_{t+1}-p_t\rangle\leq\eta_td_y\|p_{t+1}-p_t\|_*,
	\end{equation}
	where the last inequality is due to the Cauchy-Schwarz inequality and the uniform boundedness assumption in Assumption \ref{assmpt:regularity}. Strong convexity of $\Phi$ implies
	\begin{equation}\label{ineq:strong_monot}
		\langle\nabla\Phi(p_{t+1})-\nabla\Phi(p_{t}),p_{t+1}-p_t \rangle\geq\|p_{t+1}-p_t\|_*^2.
	\end{equation}
	Combining inequalities \eqref{ineq:grad_p} and \eqref{ineq:strong_monot} yields the fourth inequality.
	~\qed
\end{proof}

\section{Proofs}
\subsection{Proof of Lemma \ref{lem:Lag}}
\begin{proof}
	Note that by definition of $f_t^*$ and compactness of $\Omega_t$, \begin{displaymath}
		\dom\left(\sum_{t=1}^T f_t^*\right)=\bR^m\supseteq\dom(Th_{\Psi}),
	\end{displaymath}
	where $\dom(\cdot)$ denotes the effective domain of a function.
	By \cite[Corollary 31.2.1]{rockafellar1970convex}, it suffices to show that there exists $(\by_t)_{t=1}^T\in\ri(\proj_{\by}(\prod_{t=1}^T\conv(\Omega_t)))$ such that $\sum_{t=1}^T\by_t\in\ri(T\Psi)$, where $\ri(\cdot)$ denotes the relative interior of a set. Due to the regularity assumptions, we have\begin{align}\label{eqn:ri}
		\sum_{t=1}^T\tilde{\by}_t+\frac{T\underline{d}}{2}B_{\|\cdot\|}\subseteq \ri(T\Psi).
	\end{align}
	Choose arbitrary $(\check{\by}_t)_{t=1}^T\in\ri(\proj_{\by}(\prod_{t=1}^T\conv(\Omega_t)))$. Then, we have \begin{displaymath}
		(\by_t^\lambda)_{t=1}^T:=(1-\lambda)(\tilde{\by}_t)_{t=1}^T+\lambda(\check{\by}_t)_{t=1}^T\\
		\in\ri\Big(\proj_{\by}\Big(\prod_{t=1}^T\conv(\Omega_t)\Big)\Big)
	\end{displaymath}
	for all $\lambda\in(0,1]$ by \cite[Theorem 6.1]{rockafellar1970convex}, and $\sum_{t=1}^T\by_t^\lambda\in\ri(T\Psi)$ for some small enough $\lambda>0$ due to \eqref{eqn:ri}. The conclusion then follows. ~\qed
\end{proof}
\subsection{Proof of Lemma \ref{lem:boundedNorm}}
\begin{proof}
	Define $\kappa:=\underline d-\delta_\times=\Omega(1)$ and $b:=d_r+\delta_+=O(1)$. For all $t$, define
	$$g_t:=v_t-y_t,~D_t:=D_{\Phi}(0,p_t).$$
	Combining the second inequality with $p=0$ and the third inequality of Lemma \ref{lem:ineq}, we have
	\begin{equation}\label{pD_ineq}
		\kappa\|p_t\|_*-b
		\le
		\frac{D_t-D_{t+1}}{\eta_t}
		+
		\frac{\eta_t}{2}\|g_t\|^2\le\frac{D_t-D_{t+1}}{\eta_t}
		+
		\frac{\eta_t}{2}d_y^2,
	\end{equation}
	where the last inequality is due to Assumption \ref{assmpt:regularity}.
	
	We next consider two cases. If
	$\|p_t\|_*\geq B:=\kappa^{-1}{(b+\bar{\eta}d_y^2/2)}=O(1)$, then inequality \eqref{pD_ineq} implies
	$$D_{t+1}\leq D_t-\eta_t(\kappa \|p_t\|_*-b)+\frac{\eta_t^2}{2}d_y^2\leq D_t.$$
	Otherwise, $$\|p_{t+1}\|_*\leq\|p_t\|_*+\|p_{t+1}-p_t\|_*\leq B+\bar{\eta}d_y,$$
	where the last inequality is due to the fourth inequality of Lemma \ref{lem:ineq}, implying that
	$$D_{t+1}\leq \frac{L}{2}\|p_{t+1}\|_*^2\leq\frac{L(B+\bar{\eta}d_y)^2}{2},$$
    where the first inequality is due to Proposition \ref{prop:D_norm}.
	Combining both cases, by induction and Proposition \ref{prop:D_norm}, for all $t$, we have
	$$\|p_t\|_*\leq \sqrt{2 D_t}\leq \sqrt{L(B+\bar{\eta}d_y)^2}=O(1).~\qed$$
\end{proof}

\subsection{Proof of Theorem \ref{thm:boundedVio}}
\begin{proof}
    By the first-order optimality condition of \eqref{step:mirror_descent}, for $t=1,\ldots,T$, there exists $\xi_t\in N_{C^\circ}(p_{t+1})$ such that\begin{equation}\label{eq:bound_dist}
        y_t-v_t-\frac{\xi_t}{\eta_t}=\frac{\left(\nabla\Phi(p_{t+1})-\nabla\Phi(0)\right)-\left(\nabla\Phi(p_t)-\nabla\Phi(0)\right)}{\eta_t}.
    \end{equation}
    For all $t$, define $\Delta_t:=\nabla\Phi(p_t)-\nabla\Phi(0)$. Then, by $L$-smoothness of $\Phi$ and Lemma \ref{lem:boundedNorm}, $\|\Delta_t\|\leq L\|p_t\|_*=O(1)$ for all $t$.

    Note that $C$ is a closed convex cone. Therefore, for all $t$, we have $\xi_t\in N_{C^\circ}(p_{t+1})\subseteq C$, and
    \begin{equation}\label{v_xi:decompose}
        v_t+\frac{\xi_t}{\eta_t}\in Q+C=\Psi.
    \end{equation}
    Combining \eqref{eq:bound_dist} and \eqref{v_xi:decompose}, we have
    \begin{displaymath}
        \sum_{t=1}^T\by_t-\sum_{t=1}^T\Big(\bv_t+\frac{\xi_t}{\eta_t}\Big)=\sum_{t=1}^T\frac{\Delta_{t+1}-\Delta_t}{\eta_t}=\sum_{t=2}^T\Big(\frac{1}{\eta_{t-1}}-\frac{1}{\eta_t}\Big)\Delta_t+\frac{1}{\eta_{T}}\Delta_{T+1}.
    \end{displaymath}
    Following $\|\Delta_t\|=O(1)$ for all $t$,
    \begin{multline*} \GoalVio=\dist_{\|\cdot\|}\Big(\sum_{t=1}^T\by_t,T\Psi\Big)
        \leq\left\|\sum_{t=1}^T\by_t-\sum_{t=1}^T\Big(\bv_t+\frac{\xi_t}{\eta_t}\Big)\right\|\\
        \leq\sum_{t=2}^T\Big|\frac{1}{\eta_{t-1}}-\frac{1}{\eta_t}\Big|\cdot\|\Delta_t\|+\frac{1}{\eta_{T}}\cdot\|\Delta_{T+1}\|
        =\sum_{t=2}^TO\Big(\frac{1}{\sqrt{t}}\Big)+O(\sqrt{T})
        =O(\sqrt{T}).~\qed
    \end{multline*}
\end{proof}

\subsection{Proof of Theorem \ref{thm:main}}
\begin{proof}
Let $({r}^*_t,{x}^*_t,{y}^*_t)_{t=1}^{T}$ denote the optimal solution of \eqref{offline_convex_primal} in Lemma \ref{lem:Lag} with $(\Omega_t)_{t=1}^{T}=(Z_t)_{t=1}^{T}$. Then, under the uniform random permutation model, with $(\Omega_t)_{t=1}^{T}=(Z_{\pi(t)})_{t=1}^{T}$ and $\pi$ being a uniform random permutation of $\{1,\ldots,T\}$, we have
\begin{align*}
	\Reward= &\sum_{t=1}^Tr_t\\
    =&\sum_{t=1}^T \left[\max_{(r,x,y)\in\Omega_t}\bigl\{ r - \langle p_t, y \rangle \bigr\}+ \langle p_t, y_t \rangle -e_t\right]\\
    =&\sum_{t=1}^T \left[\max_{(r,x,y)\in\conv(\Omega_t)}\bigl\{ r - \langle p_t, y \rangle \bigr\}+ \langle p_t, y_t \rangle\right]-\sum_{t=1}^Te_t\\
    \geq & ~ \sum_{t=1}^{T}{r}^*_{\pi(t)}+\sum_{t=1}^{T}\langle\bp_t, \by_t-{\by}^*_{\pi(t)}\rangle-\sum_{t=1}^Te_t\\
	= & ~ z_R+\sum_{t=1}^{T}\langle \bp_t,\by_t-\bv_t\rangle+\sum_{t=1}^{T}\langle\bp_t,\bv_t-\bEpiO[{\by}^*_{\pi(t)}]\rangle~+\\
	&\quad\sum_{t=1}^{T}\langle\bp_t,\bEpiO[{\by}^*_{\pi(t)}]-{\by}^*_{\pi(t)}\rangle-\sum_{t=1}^Te_t.
\end{align*}
Here, the second equality is due to the definition \eqref{define_et} of $e_t$, and the first inequality holds since $({r}^*_{\pi(t)},{\bx}^*_{\pi(t)},{\by}^*_{\pi(t)})\in\conv(\Omega_t)$ from Lemma~\ref{lem:Lag}.

To bound $\bEpiO[\Reward]$, we bound the expectations of $\sum_{t=1}^{T}\langle \bp_t,\by_t-\bv_t\rangle$, $\sum_{t=1}^{T}\langle\bp_t,\bv_t-\bE_\pi[{\by}^*_{\pi(t)}]\rangle$ and $\sum_{t=1}^T e_t$ separately. By rearranging the second inequality of Lemma \ref{lem:ineq} with $\bp=\mathbf{0}$, we have\begin{align*}
	\langle \bp_t,\by_t-\bv_t\rangle \geq \frac{D_\Phi(0,p_{t+1})-D_\Phi(0,p_t)}{\eta_t}-\frac{\eta_t}{2}{d_y^2}.
\end{align*}
Note that $\bEpiO[D_\Phi(0,p_t)]=O(1)$ for all $t$, following Proposition \ref{prop:D_norm} and Lemma \ref{lem:boundedNorm}. It then follows that
\begin{align*}
	\bEpiO\left[\sum_{t=1}^{T}\langle \bp_t,\by_t-\bv_t\rangle\right]
	\geq&\left(\sum_{t=1}^T \frac{[\bEpiO[D_\Phi(0,p_{t+1})]-\bEpiO[D_\Phi(0,p_t)]}{\eta_t}-\frac{\eta_t}{2}{d_y^2}\right)\\
    &\geq-\sum_{t=2}^T(\frac{1}{\eta_t}-\frac{1}{\eta_{t-1}})\bEpiO[D_\Phi(0,p_{t})]-\sum_{t=1}^{T}\frac{\eta_t}{2}{d_y^2}=-O(\sqrt{T}).
\end{align*}
Also note that $\bEpiO[{\by}^*_{\pi(t)}] =   T^{-1}\sum_{\tau=1}^{T}{\by}^*_\tau \in \Psi$
for all $t$ as $\pi$ is a uniform random permutation. It implies that $\langle\bp_t,\bv_t-\bEpiO[{\by}^*_{\pi(t)}]\rangle=h_{{\Psi}}(\bp_t)-\langle\bp_t,\bEpiO[{\by}^*_{\pi(t)}]\rangle \geq 0$ for all $t$ by the definition of $h_{{\Psi}}(\cdot)$. Moreover, following Assumption \ref{assmpt:approx} and Lemma \ref{lem:boundedNorm}, $\bEpiO[e_t]\leq \delta_++\delta_\times\bEpiO[\|p_t\|_*]\leq \delta_++O(\delta_\times)$.
Therefore,\begin{equation}
	\bEpiO[\Reward]\geq~z_R - O(\sqrt{T})-\delta_+ T-O(\delta_\times T) +\sum_{t=1}^{T}\bEpiO[\langle\bp_t,\bE_\pi[{\by}^*_{\pi(t)}]-{\by}^*_{\pi(t)}\rangle].\label{RewardBnd:deterministic}
\end{equation}
Let $\mF_{t-1}$ denote the sigma algebra generated by the random events up to time step $t-1$.
Next, we bound the expectation of $\langle\bp_t,\bE_\pi[{\by}^*_{\pi(t)}]-{\by}^*_{\pi(t)}\rangle$ conditional on  $\mF_{t-1}$, 
\begin{equation}\label{ineq:cond_exp}
\begin{aligned}
    \bEpiO&[\langle\bp_t,\bEpiO[{\by}^*_{\pi(t)}]-{\by}^*_{\pi(t)}\rangle|\mF_{t-1}]
	=~ \langle\bp_t,\bEpiO[{\by}^*_{\pi(t)}]-\bEpiO[{\by}^*_{\pi(t)}|\mF_{t-1}]\rangle\\
	=&~ \left\langle\bp_t,\bEpiO[{\by}^*_{\pi(t)}]-\frac{1}{T-t+1}\sum_{\tau=t}^{T}{\by}^*_{\pi(\tau)}\right\rangle
	\geq~ -\|\bp_t\|_*\left\|\bEpiO[{\by}^*_{\pi(t)}]-\frac{1}{T-t+1}\sum_{\tau=t}^{T}{\by}^*_{\pi(\tau)}\right\|,
\end{aligned}
\end{equation}
where the first equality holds since $\bp_t$ is determined based on $\mathcal{F}_{t - 1}$, the second equality is due to the fact that $\pi$ is a uniform random permutation.

Let $\alpha, \beta>0$ be $\|\cdot\|$-dependent constants such that
$$\|\cdot\|_\infty\leq \alpha\|\cdot\|,~\|\cdot\|\leq\beta\|\cdot\|_\infty.$$ Fix an arbitrary $v_0\in Q$.
Then the uniform boundedness assumption in Assumption \ref{assmpt:regularity} implies $\|y-v_0\|_\infty\leq \alpha\|y-v_0\|\leq \alpha d_y.$
	As $\pi$ is uniformly chosen at random, $\big(\sigma(t)=\pi(T-t+1)\big)_{t=1}^{T}$ is also a uniform random permutation of $\mT$.   
	By Hoeffding's inequality for sampling without replacement \cite{10.2307/2282952}, for all $\epsilon>0$ and $i\in\{1,\ldots,m\}$, we have\begin{displaymath}
		\bP_{\pi}\left(\left|\left(\bEpiO[{\by}^*_{\pi(t)}]-\frac{1}{T-t+1}\sum_{\tau=t}^{T}{\by}^*_{\pi(\tau)}\right)_i\right|>\epsilon\right)
		\leq2\exp\Big(-\frac{(T-t+1)\epsilon^2}{2{\alpha^2 d_y^2}}\Big).
	\end{displaymath}
	In other words, for all $\rho\in(0,1]$, with probability at least $1-\rho/m$,\begin{displaymath}
		\left|\left(\bEpiO[{\by}^*_{\pi(t)}]-\frac{1}{T-t+1}\sum_{\tau=t}^{T}{\by}^*_{\pi(\tau)}\right)_i\right|
        \leq \sqrt{\frac{2\alpha^2 {d_y^2}\log(2m/\rho)}{T-t+1}}.
	\end{displaymath}
	Then, by taking the union bound over $i\in\{1,\ldots,m\}$, with probability at least $1-\rho$,\begin{multline*}
		\left\|\bEpiO[{\by}^*_{\pi(t)}]-\frac{1}{T-t+1}\sum_{\tau=t}^{T}{\by}^*_{\pi(\tau)}\right\|\leq\beta\left\|\bEpiO[{\by}^*_{\pi(t)}]-\frac{1}{T-t+1}\sum_{\tau=t}^{T}{\by}^*_{\pi(\tau)}\right\|_\infty
        \\\leq\alpha\beta d_y\cdot\sqrt{\frac{2\log(2m/\rho)}{T-t+1}}
        =O_{\|\cdot\|,m}\Big(\sqrt{\frac{\log(2/\rho)}{T-t+1}}\Big).
	\end{multline*}
	By inequality \eqref{ineq:cond_exp}, Lemma \ref{lem:boundedNorm} and integrating the quantile function, we have\begin{align}
		\bEpiO[\langle\bp_t,\bEpiO[{\by}^*_{\pi(t)}]&-{\by}^*_{\pi(t)}\rangle]
        =~\bEpiO\left[\bEpiO[\langle\bp_t,\bE_\pi[{\by}^*_{\pi(t)}]-{\by}^*_{\pi(t)}\rangle|\mF_{t-1}]\right]\nonumber\\
		\geq&~-\bEpiO\left[\|\bp_t\|_*\left\|\bEpiO[{\by}^*_{\pi(t)}]-\frac{1}{T-t+1}\sum_{\tau=t}^{T}{\by}^*_{\pi(\tau)}\right\|\right]\nonumber\\
		&\geq~-\int_{0}^1 O_{\|\cdot\|,m}\Big(\sqrt{\frac{\log(2/\rho)}{T-t+1}}\Big)~d\rho
        =~- O_{\|\cdot\|,m}\Big(\sqrt{\frac{1}{T-t+1}}\Big),\label{ExpBnd}
	\end{align}
    where the second inequality uses the fact that Lemma \ref{lem:boundedNorm} holds for any arrival order.
It then follows from \eqref{RewardBnd:deterministic} that\begin{multline*}
    \bEpiO[\Reward]
        \geq z_R+\sum_{t=1}^{T}\bEpiO[\langle\bp_t,\bE_\pi[{\by}^*_{\pi(t)}]-{\by}^*_{\pi(t)}\rangle]-O(\sqrt{T})-\delta_+ T-O(\delta_\times T)\\
		\geq z_R-O_{\|\cdot\|,m}(\sqrt{T})-\delta_+ T-O(\delta_\times T).
~\qed
	\end{multline*}
\end{proof}

\section{Numerical Details} \label{apd:numerical}

\subsection{Random Instance Generation}\label{apd:details-random-instances}
In the MNL model, the reward function (i.e., expected revenue) of customer $t$ is given by
\[
r_t(S)=\frac{\sum_{i\in S}\gamma_i\nu_{ti}}{1+\sum_{i\in S}\nu_{ti}}.
\]
For each base instance, we independently sample $\gamma_i\sim\mathrm{Unif}[0.7,1]$ for $i\in\{1,2,3,4\}$ and $\gamma_i\sim\mathrm{Unif}[0.1,0.3]$ for $i\in\{5,6,7,8\}$. Each $f_i$ is sampled uniformly from $\{0.20,0.22,0.24,0.26,0.28,0.30\}$, every $\nu_{ti}$ is independently sampled from $\mathrm{Unif}[0.2,2]$.

\subsection{Results for A Deliberately Constructed Instance} \label{apd:artificial-instance}

This constructed instance examines the two \GoalVio\ regimes given in Theorem~\ref{thm:boundedVio} and Example~\ref{exmp:goalvio} using the FPTAS proposed by \cite{cao2023revenue}. It has $52$ products, with at most $s=51$ displayed per customer. The individual-product requirements in Section~\ref{sec:numerical} are replaced by one aggregate display requirement for products $1,\ldots,51$:
\begin{equation}\label{eq:num-artificial}
 \max_{(r_t, x_t, y_t) \in \Omega_t ~ \forall t \in [T]}\ \sum_{t=1}^T r_t
 \quad\text{s.t.}\quad \sum_{t = 1}^T y_t \geq 50.1 T 
\end{equation} 
with local constraint set 
\begin{align*}
    \Omega_t := \left\{ (r, x, y) ~:~ \sum_{i} x_i \leq 51, ~ r = \frac{\sum_i \gamma_i \nu_{ti} x_i }{1 + \sum_i \nu_{ti} x_i }, ~ y = \sum_{i = 1}^{51} x_i, ~ x \in \{0, 1\}^{52} \right\}. 
\end{align*}
Thus $y_t$ counts target-group $\{1,\ldots,51\}$ displays, and product $52$ does not contribute. The no-purchase attraction is $1$. All product coefficients are given in Table~\ref{tab:num-artificial-catalog}. 


\begin{table}[H]\centering
\caption{Coefficients and types of customers for \eqref{eq:num-artificial}.}
\label{tab:num-artificial-catalog}
\begin{tabular}{@{}lrrrc@{}}\toprule
Products & ~~Revenue $\gamma_i$ & ~~Low attraction & ~~High attraction & ~~Target group?\\\midrule
$1,\ldots,50$ & $1$ & $0.02$ & $0.04$ & Yes\\
$51$ & $0.1$ & $0.2$ & $0.4$ & Yes\\
$52$ & $0.1$ & $0.002$ & $0.004$ & No\\\bottomrule
\end{tabular}
\end{table}

We consider an adversarial arriving order as Theorem~\ref{thm:boundedVio} does not depend on the arriving order. For each $T\in\{500,2000,8000\}$, the first $T/2$ customers have low attraction and the rest have high attraction. 

\noindent
\textbf{Algorithm and oracle settings.}
With $\Psi=[50.1,\infty)$ and norm $|\cdot|$, the maximum group exposure $51$ gives slack $\underline d=51-50.1=0.9$. We use $\Phi(p)=p^2/2$, $p_1=0$, and $\eta_t=0.2/\sqrt t$. Since $C^\circ=\mathbb R_-$, the updating rule in Algorithm~\ref{alg:LagDual} becomes
\[
 (r_t,x_t,y_t)=\mathcal O_t(p_t;\Omega_t),\qquad
 p_{t+1}=-[-p_t+\eta_t(50.1-y_t)]_+.
\]
Each period calls the FPTAS in \cite[Theorem~5.2 and Appendix~F]{cao2023revenue} to approximately maximize $r-p_ty$.

We fix $\delta_+=0.03$ and consider two settings of the
FPTAS in \cite{cao2023revenue}. We consider relative approximation errors $\epsilon=0.8888$ and $\epsilon=0.9$ in the FPTAS, corresponding to $\alpha=0.9$ (i.e., $\delta_\times=0.81$) and $\alpha=10/9$ (i.e., $\delta_\times=1$), respectively.

\begin{figure}[h!]
	\centering
	{\includegraphics[width=0.48\textwidth]{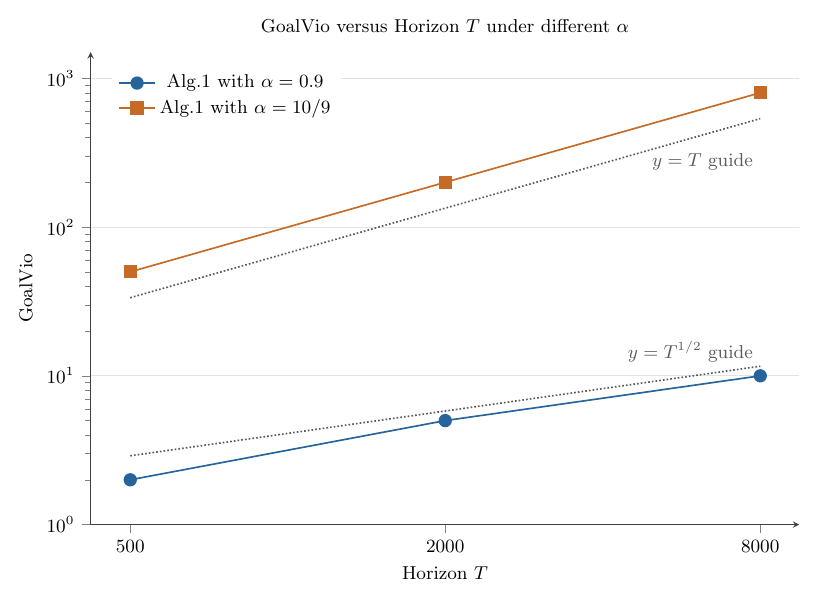}}
    {\includegraphics[width=0.48\textwidth]{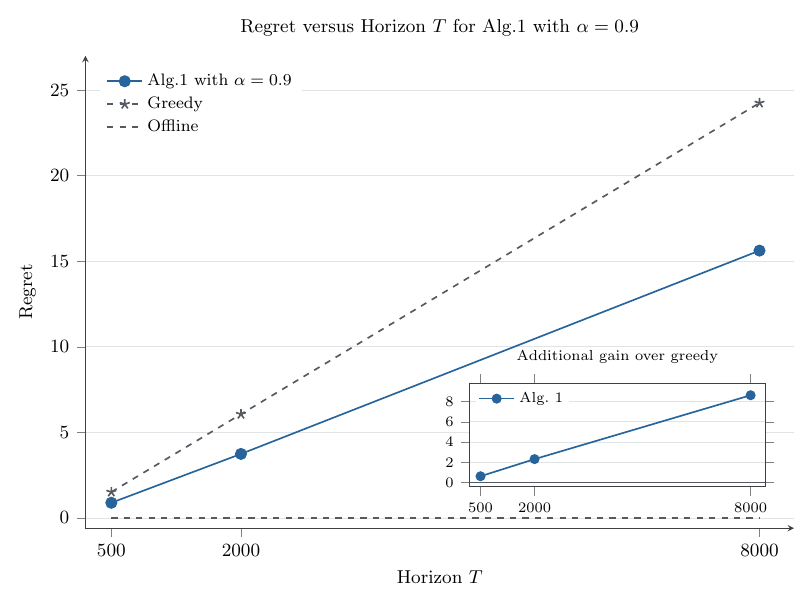}}

    \vspace{- 3 ex} 
    \caption{\small This figure illustrates the dependencies of GoalVio and Regret over the deliberately constructed instance. Left panel: GoalVio versus horizon $T$. Right panel: Regret and additional gain over Greedy versus horizon $T$.}
	\label{fig:constructed-instance}
l\end{figure}




\noindent
\textbf{Numerical results and discussions.} 
Figure~\ref{fig:constructed-instance} reports the dependencies of \GoalVio\ and \Reward\  (regret) over the deliberately constructed instance. 

In the left panel, \GoalVio\ grows approximately as $\sqrt{T}$ at $\alpha=0.9$ and linearly at $\alpha=10/9$, consistent with the phase-transition behavior illustrated in Example~\ref{exmp:goalvio}. Here the subproblems are solved using the FPTAS proposed in \cite{cao2023revenue}. Under our specified mapping from $\alpha$ to oracle precision, the transition occurs between $1$ and $10/9$. Thus, the sufficient condition $\delta_{\times}<\underline{d}$ need not identify the exact transition point for a particular instance and oracle.

We plot \Reward\  (regret) only for $\alpha=0.9$ as the $\alpha=10/9$ solution violates the visibility constraint too much. The right panel shows that Algorithm~\ref{alg:LagDual} achieves lower regret, and hence higher reward, than Greedy at every tested horizon. 



\end{document}